\documentclass[dvips,ejs]{imsart}

\RequirePackage[OT1]{fontenc}
\RequirePackage{amsthm,amsmath,natbib}
\RequirePackage[colorlinks,citecolor=blue,urlcolor=blue]{hyperref}
\RequirePackage{hypernat}

\usepackage{mathrsfs}
\usepackage{amsfonts}
\usepackage{amssymb}
\usepackage{amsmath,epsfig}

\pubyear{2009}
\startlocaldefs
\numberwithin{equation}{section}
\theoremstyle{plain}

\endlocaldefs

\startlocaldefs
\numberwithin{equation}{section}
\theoremstyle{plain}

\newtheorem{lemma}{Lemma}[section]

\endlocaldefs

\begin{document}

\begin{frontmatter}
\title{A Bernstein-type inequality for stochastic processes
of quadratic forms of Gaussian variables
%\thanksref{T1}
}
\runtitle{A Bernstein-type inequality for quadratic forms of gaussian variables}
%\thankstext{T1}{Footnote to the title with the `thankstext' command.}

\begin{aug}
\author{\fnms{Ikhlef} \snm{Bechar}\thanksref{t2}
%\ead[label=e1]{ikhlef.bechar@cardiov.ox.ac.uk}
\ead[label=e2]{ikhlef.bechar@sophia.inria.fr} }
%\and
%\author{\fnms{Second} \snm{Author}\thanksref{t3}\ead[label=e2]{second@somewhere.com}}

%\address{Address of the First and Second authors\\
%usually few lines long\\
%\printead{e1,e2}}

%\author{\fnms{Third} \snm{Author}
%\ead[label=e3]{third@somewhere.com}
%\ead[label=u1,url]{www.foo.com}}

%\address{Address of the Third author\\
%usually few lines long\\
%usually few lines long\\
%\printead{e3}\\
%\printead{u1} }

%\thankstext{t1}{EPSRC Postdoctoral Research Grant, Oxford University}
\thankstext{t2}{Pulsar project, INRIA Sophia Antiplois}
%\thankstext{t3}{Second supporter of the project}
\runauthor{I. Bechar}

%\affiliation{Some University and Another University}

\runauthor{I. Bechar}

\affiliation{INRIA Nice Sophia Antipolis}

%\address{ OCMR / FMRIB  \\
%Oxford University\\
%John Radcliffe Hospital\\
%Oxford, OX3 9DU \\
%United Kingdom \\
%\printead{e1}\\
%\phantom{E-mail: ikhlef.bechar@cardiov.ox.ac.uk\ } }

\address{Projet Pulsar \\
 INRIA Sophia Antipolis \\
Route des Lucioles - BP 93 \\
06902 Sophia Antipolis Cedex \\
FRANCE \\
\printead{e2}\\
\phantom{E-mail: ikhlef.bechar@sophia.inria.fr \ } }
\end{aug}

\begin{abstract}
%\selectlanguage{english}
% Text of abstract in English
We introduce a Bernstein-type inequality which serves to uniformly control quadratic forms of
gaussian variables. The latter can for example be used to derive sharp model selection 
criteria for linear estimation in linear regression and linear inverse problems via 
penalization, and we do not exclude that its scope of application can be made even broader.
\end{abstract}
\end{frontmatter}

% now the Version française abrégée, if it exists
%\selectlanguage{francais}

%\selectlanguage{english}
% main text
\section*{A Bernstein-type inequality for quadratic forms of gaussian variables}
The concentration phenomenon of stochastic processes around their mean is of
key importance in statistical estimation by model selection for getting 
non-asymptotic bounds for some statistics. For example in model selection
via penalization, for devising sharp penalties and proving useful upper bounds
for the risk of an estimator, one needs generally to control uniformly the 
statistic of the risk of an estimator by means of a sharp concentration inequality. 
This topic has received since recently (late nineties) a considerable interest 
among the statistical community above all further to the amazing series of works of Michel Talagrand
which can be seen as the infinite dimensional analogue of the Bernstein's 
inequality (see in particular \cite{talagrand_1995} for an overview and \cite{talagrand_1996} for later 
advances). Their  application in non-asymptotic model selection has first been discovered by 
Birg\'e and Massart (e.g. \cite{birge_massart_1998}), then refined and popularized by the same authors (e.g. \cite{birge_massart_2001,birge_massart_2007}). For beautiful lectures on the topic, we refer
the dear reader to \cite{massart_2000, massart_2007}.

In this small body of work, we establish a new Berstein-type inequality which serves
to control (e.g. uniformly) quadratic forms of Gaussian variables and which happens 
to be useful for controlling, for example, uniformly the quadratic risk of 
a finite (or a countable) set of linear estimators in linear regression and linear 
inverse problems (see \cite{bechar_2009a} for an application). In the remainder, 
we will give both the uncorrelated form and the correlated form of such an inequality.
{\begin{lemma}
\label{lemma_1}
Let $a =(a_k)_{k=1,p}$  and $b=(b_k)_{k=1,p}$ be two $p-$dimensional real
vectors, and  consider the following random expression : $T = \sum_{k=1}^p a_k z_k^2 + b_k z_k$,
where $z_k, k=1,\cdots,p$ are i.i.d. $N(0,1)$, and let's put : $a^{+} = \sup\{\sup_{k=1,\cdots,p}\{a_k\}, 0\}$,
$a^{-} = \sup\{\sup_{k=1,\cdots,p}\{-a_k\}, 0\}$. Then the following two concentration results
hold true for all $x>0$ :
\begin{equation}
\label{conc1}
\mathbb{P}\Big [T \geq \sum_{k=1}^p a_k  +  2 \sqrt{\sum_{k=1}^p a_k^2 + \frac{b_k^2}{2}} \sqrt{x} + 2 a^{+}x \Big] \leq \exp[-x]
\end{equation}
\begin{equation}
\label{conc2}
\mathbb{P}\Big [T \leq \sum_{k=1}^p a_k  -  2 \sqrt{\sum_{k=1}^p a_k^2 + \frac{b_k^2}{2}} \sqrt{x} - 2 a^{-}x \Big] \leq \exp[-x]
\end{equation}
\end{lemma}
The proof of the lemma is rather technical, so it is deferred to the appendix section.}\\
The following lemma uses the concentration results of lemma (\ref{lemma_1})
to control more general quadratic forms of Gaussian variables involving a matrix.
\begin{lemma}
\label{theo_1}
Consider the random expression $T =  z^T A z + b^T z$, where $A$ is $p$ by $p$ real square matrix,
$b$ is a $p-$dimensional real vector, and $z = (z_k)_{k=1,p}$ is a $p-$dimensional standard gaussian vector,
i.e. $z_k, k=1,p$ are i.i.d. zero-mean gaussian variables with standard deviation $1$. Let's denote by $s_k, k=1,p$
the eigen values of the symmetric matrix $\frac{1}{2}\big(A + A^T\big)$, and let's put $s^+ = \sup\{\sup_{k=1,\cdots,p}\{s_k\}, 0\}$,
and $s^{-} = \sup\{\sup_{k=1,\cdots,p}\{-s_k\}, 0\}$. Then, the following two concentration results
hold true for all $x>0$
\begin{equation}
\label{conc1b}
\mathbb{P}\Big [T \geq tr(A)  +  2 \sqrt{\frac{1}{4}\|A + A^T\|^2 + \frac{1}{2}\|b\|^{2}} \sqrt{x} + 2 s^{+}x \Big] \leq \exp[-x]
\end{equation}
\begin{equation}
\label{conc2b}
\mathbb{P}\Big [T \leq tr(A)  -  2 \sqrt{ \frac{1}{4}\|A + A^T\|^2 + \frac{1}{2}\|b\|^{2}} \sqrt{x} - 2 s^{-}x \Big] \leq \exp[-x]
\end{equation}
\end{lemma}
{\it Proof.} {
One can rewrite $T$ as follows: $T =  z^T A z + b^T z =  \frac{1}{2} z^T \big (A+A^T\big) z + b^T z$, and by using the eigen value
decomposition of the symmetric matrix $\frac{1}{2}\big(A+A^T\big)$ one derives $T =  \sum_{k=1}^p s_k {{z'}_k^2} + {b'}_k z'$,
where $s_k, k=1,p$ are the respective eigen values of $\frac{1}{2}\big(A+A^T\big)$, $z' = U^T z$  with $U$ standing for the
(orthonormal) eigen matrix of $\frac{1}{2}\big(A+A^T\big)$,  and $b' = U^T b$. Then, by noticing that $z'$ stands for a
$p-$dimensional standard gaussian vector, $\|b'\|^2 = \|b\|^2$, $\sum_{k=1}^p s_k = tr(A)$,
and $\sum_{k=1}^p s_k^2 =  \frac{1}{4} \big\|A+A^T \big\|^2$, so by applying lemma (\ref{lemma_1}), the proof
of lemma (\ref{theo_1}) follows immediately.
%\end{proof}
}

% etc, etc

 \appendix
 \label{app1}
\section{Proof of Lemma (\ref{lemma_1})}
\begin{proof}
We make use of the following lemma for proving lemma (\ref{lemma_1})
\begin{lemma}[Birge \& Massart 1998]
\label{lemma_1c}
If a random variable $\xi$ satisfies for some two real positive numbers $u$ and $v$ the following inequality :
\begin{equation}
\log\Big( \mathbb{E}\Big[\exp[ y \xi ] \Big] \Big) \leq \frac{(u y)^2}{1-v y}, \textrm{for all \,\,\,} 0<y<\frac{1}{v}
\end{equation}
then
\begin{equation}
\mathbb{P}\Big[ \xi \geq 2 u\sqrt{x} + v x \Big] \leq \exp[-x], \textrm{for all \,\,\,} x>0
\end{equation}
\end{lemma}
We refer the reader to \cite{birge_massart_1998} for a proof of this lemma. \\

Now, to prove lemma (\ref{lemma_1}), one can notice first that concentration inequality (\ref{conc2}) can
be obtained from (\ref{conc1}) by considering the random quantity
\begin{equation}
T' = -T = \sum_{k=1}^p (-a_k) z_k^2 + (-b_k) z_k \nonumber
\end{equation}
and by applying (\ref{conc1b}) on $T'$ instead of $T$. So, we need to prove only (\ref{conc1b}).
To do this, let us rewrite $T$ as follows: $ T= \sum_{k=1}^p T_k$, where $ T_k =  a_k z_k ^2 + b_k z_k$,
and let us compute $\log\big [\mathbb{E}\big (\exp(y (T - \bar{T}))) \big] $, where $\bar{T} = \sum_{k=1}^p a_k$.
We have
\begin{displaymath}
\mathbb{E}\big[\exp(y T_k)\big] = \frac{1}{\sqrt{2\pi}}\int_{-\infty}^{\infty} \exp\Big [-\frac{1}{2} ((1-2 a_k y) t^2 - 2 y b_k t) \Big] dt
\end{displaymath}

\begin{displaymath}
\nonumber
\mathbb{E}\big[\exp[y T_k]\big] = \exp\bigg[ \frac{ \frac{b_k^2}{2} y^2 }{1-2 a_k y}  \bigg] \bigg (\frac{1}{\sqrt{2\pi}}\int_{-\infty}^{\infty} \exp\bigg [-\frac{1}{2} \Big ( \sqrt{1-2 a_k y} t - \frac{ b_k y}{\sqrt{1-2 a_k y}}\Big )^2 \bigg] \mathrm{d} t \bigg )
\end{displaymath}

\begin{displaymath}
\mathbb{E}\big[\exp[y T_k]\big] = \frac{\exp\bigg[ \frac{\frac{b_k^2}{2} y^2}{1-2 a_k y} \bigg] }{\sqrt{1-2 a_k y}}
\end{displaymath}

\begin{displaymath}
\mathbb{E}\Big[\exp[y (T_k - a_k )]\Big] = \frac{\exp\bigg[ \frac{\frac{b_k^2}{2} y^2 }{1-2 a_k y} \bigg] \exp[-y a_k] }{\sqrt{1-2 a_k y}}
\end{displaymath}

\begin{displaymath}
\log \Big (\mathbb{E}\Big[\exp[y (T_k - a_k )]\Big] \Big) =  \frac{\frac{b_k^2}{2} y^2}{1-2 a_k y} - \frac{1}{2}\log \Big(1-2 a_k y\Big) - a_k y
\end{displaymath}
Then, by putting $a^{+} = \sup\big\{ \sup_{k=1,\cdots,p}\{a_k\},0 \big\}$, one derives (see the technical details below) that for all
 $0<y<\frac{1}{2a^{+}}$
\begin{displaymath}
\log \Big (\mathbb{E}\Big[\exp[y (T_k - a_k )]\Big] \Big) \leq \frac{\big (a_k^2 + \frac{b_k^2}{2}\big) y^2}{1-2 a_{k} y} \leq \frac{\big (a_k^2 + \frac{b_k^2}{2}\big) y^2}{1-2 a^{+} y}
\end{displaymath}
which implies by independence that for all $0<y<\frac{1}{2a^{+}}$
\begin{displaymath}
\log \Big (\mathbb{E}\Big[\exp[y (T - \bar{T} )]\Big] \Big) \leq \sum_{k=1}^{p}  \frac{\big (a_k^2 + \frac{b_k^2}{2}\big) y^2}{1-2 a^{+} y}
\end{displaymath}

\begin{displaymath}
\log \Big (\mathbb{E}\Big[\exp[y (T - \bar{T} )]\Big] \Big) \leq   \frac{\Big (\sum_{k=1}^{p} \big(a_k^2 + \frac{b_k^2}{2}\big) \Big) y^2}{1-2 a^{+} y}
\end{displaymath}
 Finally, by applying lemma (\ref{lemma_1c}) below with $ u = \sqrt{\sum_{k=1}^{p} (a_k^2 + \frac{b_k^2}{2})} $ , and $ v = 2 a^{+}$ , one derives that for all $x>0$ :
\begin{displaymath}
\mathbb{P}\Big [ T \geq  \big [\sum_{k=1}^{p} a_k \big ] +  2 \sqrt{\sum_{k=1}^{p} \big(a_k^2 + \frac{b_k^2}{2}\big)} \sqrt{x} + 2 a^{+} x \Big ]  \leq \exp[-x]
\end{displaymath}
This terminates the proof of lemma (\ref{lemma_1})
\end{proof}

\subsection*{Some additional technical details about the proof}
\label{tech_det}
We will show here that for all $r>0$, $a\geq r$, and $0<y<\frac{1}{2a}$,
one has
\begin{equation}
\label{ineq_1}
 \frac{-1}{2} \log(1 - 2ry) - ry \leq \frac{r^2 y^2}{1-2{a} y}
\end{equation}
and that for all $r\leq 0$, for all $a>0$, and for all $0<y<\frac{1}{2a}$, one has
\begin{equation}
\label{ineq_2}
 -\frac{1}{2} \log(1 - 2ry) - ry \leq \frac{r^2 y^2}{1-2{a} y}
\end{equation}
\begin{proof}
let us start by showing inequality (\ref{ineq_1}). To do this, let us consider
the following function
\begin{equation}
f_{r, {a}}(y) = -\frac{1}{2} \log(1 - 2ry) - ry - \frac{r^2 y^2}{1-2{a} y} \nonumber
\end{equation}
One first notices that $f_{r, {a}}(0) =0$, then a sufficient condition for inequality
(\ref{ineq_1}) to hold true is  that $f_{r, {a}}(y)' \leq 0$, for all $0<y<\frac{1}{2a}$.
We have
\begin{displaymath}
f_{r, {a}}(y) = \frac{-1}{2} \log(1 - 2ry) - r y + \frac{r^2y}{2{a}} + \frac{r^2}{(2{a})^2} -  \frac{\frac{r^2}{(2{a})^2}}{1-2{a} y}
\end{displaymath}
one then derives that
\begin{displaymath}
f_{r, {a}}(y)' =  \frac{r}{1 - 2ry} - r +  \frac{r^2}{2{a}} -  \frac{\frac{r^2}{2 {a}}}{(1-2{a} y)^2}
\end{displaymath}
\begin{displaymath}
f_{r, {a}}(y)' =  \frac{2 r^2 y }{1 - 2ry}  -   \frac{r^2 y }{(1-2{a} y)} - \frac{r^2 y }{(1-2{a} y)^2}
\end{displaymath}
\begin{displaymath}
f_{r, {a}}(y)' \leq  \frac{2 r^2 y }{1 - 2ry}  -   \frac{2 r^2 y }{(1-2{a} y)}
\end{displaymath}
and finally since $ \frac{1 }{(1-2{r} y)}  \leq   \frac{1}{(1-2{a} y)}$, one deduces that
\begin{displaymath}
f_{r, {a}}(y)' \leq  \frac{2 r^2 y }{(1-2{a} y)}  -   \frac{2 r^2 y }{(1-2{a} y)} = 0
\end{displaymath}
then we have shown (\ref{ineq_1}).\\
We proceed in the same way as for showing inequality (\ref{ineq_1}) to show inequality (\ref{ineq_2}).
So let us consider the following function
\begin{equation}
g_{r, {a}}(y) = \frac{-1}{2} \log(1 - 2ry) - ry - \frac{r^2 y^2}{1-2{a} y} \nonumber
\end{equation}
One first notices that $g_{r, {a}}(0) =0$, then a sufficient condition for inequality (\ref{ineq_2}) to hold
true is to that $g_{r, {a}}(y)' \leq 0$ for all $0<y<\frac{1}{2a}$. One derives that
\begin{displaymath}
g_{r, {a}}(y)' =  \frac{r}{1 - 2ry} - r +  \frac{r^2}{2{a}} -  \frac{\frac{r^2}{2 {a}}}{(1-2{a} y)^2}
\end{displaymath}

\begin{displaymath}
g_{r, {a}}(y)' =  \frac{2 r^2 y }{1 - 2ry}  -   \frac{r^2 y }{(1-2{a} y)} - \frac{r^2 y }{(1-2{a} y)^2}
\end{displaymath}
\begin{displaymath}
g_{r, {a}}(y)' \leq  \frac{2 r^2 y }{1 - 2ry}  -   \frac{2 r^2 y }{(1-2{a} y)}
\end{displaymath}
and finally, since $\frac{1}{1 - 2ry} \leq  \frac{1}{(1-2{a} y)} $, one finds that
\begin{displaymath}
g_{r, {a}}(y)' \leq  \frac{2 r^2 y }{(1-2{a} y)}  -   \frac{2 r^2 y }{(1-2{a} y)} = 0
\end{displaymath}
\end{proof}

\subsection*{Birge's \& Massart concentration inequality}
\begin{lemma}
\label{lemma_1c}
If a random variable $\xi$ satisfies for some two real positive numbers $u$ and $v$ the following inequality :
\begin{equation}
\log\Big( \mathbb{E}\Big[\exp[ y \xi ] \Big] \Big) \leq \frac{(u y)^2}{1-v y}, \textrm{for all \,\,\,} 0<y<\frac{1}{v}
\end{equation}
then
\begin{equation}
\mathbb{P}\Big[ \xi \geq 2 u\sqrt{x} + v x \Big] \leq \exp[-x], \textrm{for all \,\,\,} x>0
\end{equation}
\end{lemma}
The proof of this lemma can be found in  \cite{birge_massart_1998}.

\end{document}